\documentclass[12pt]{article}
\usepackage{amsfonts}
\usepackage{amssymb}
\usepackage{amsmath}
\usepackage{endnotes}
\usepackage[doublespacing]{setspace}

\input{tcilatex}
\let\footnote=\endnote

\begin{document}

\title{A technical note on divergence of the Wald statistic\thanks{%
This work was supported by the Willam Dow Chair in Political Economy (McGill
University), the Bank of Canada Research Fellowship, The Toulouse School of
Economics Pierre-de-Fermat Chair of Excellence, A Guggenheim Fellowship,
Conrad-Adenauer Fellowship from Alexander-von-Humboldt Foundation, the
Canadian Network of Centres of Excellence program on Mathematics of
Information Technology and Complex Systems, the Natural Sciences and
Engineering Research Council of Canada, the Social Sciences and Humanities
Research Council of Canada and the Fonds de recherche sur la soci\'{e}t\'{e}
et la culture (Qu\'{e}bec). The authors also thank the research centres
CIREQ and CIRANO for providing support and meeting space for the joint work.
We thank Purevdorj Tuvaandorj for very useful comments.}}
\author{Jean-Marie Dufour\thanks{%
William Dow Professor of Economics, McGill University, Centre
interuniversitaire de recherche en analyse des organisations (CIRANO) and
Centre interuniversitaire de recherche en \'{e}conomie quatative (CIREQ).}
\and Eric Renault \thanks{%
Brown University} \and Victoria Zinde-Walsh\thanks{%
McGill University and CIREQ}}
\maketitle
\date{}

\begin{abstract}
The Wald test statistic has been shown to diverge (Dufour et al, 2013, 2017)
under some conditions. This note links the divergence to eigenvalues of a
polynomial matrix and establishes the divergence rate.
\end{abstract}

\section{\protect\bigskip The set-up and an example of divergence}

Suppose that a $p\times 1$ parameter of interest $\bar{\theta}$ satisfies

\begin{equation*}
H_{0}:g\left( \theta \right) =0,
\end{equation*}%
where $g\left( \theta \right) $ is a $q\times 1$ vector of differentiable
functions; $g\left( \theta \right) =\left( g_{1}\left( \theta \right)
,...,g_{q}\left( \theta \right) \right) ^{\prime };$ $q\leq p.$

Let $V$ be a symmetric positive definite matrix.

\textbf{Assumption 1.} \emph{In some open set }$\Theta \subset R^{p}$\emph{\
there is a random sequence }$\hat{\theta}_{T}\in \Theta $\emph{\ and }$%
p\times p$\emph{\ random matrix sequence, }$\hat{V}_{T},$\emph{\ such that
as }$T\rightarrow \infty $\emph{\ \ }%
\begin{eqnarray*}
\sqrt{T}V^{-\frac{1}{2}}\left( \hat{\theta}_{T}-\bar{\theta}\right)
&\rightarrow &_{d}Z; \\
Z &\sim &N\left( 0,I_{p}\right) ; \\
\text{\textit{and} }\hat{V}_{T} &\rightarrow &_{p}V.
\end{eqnarray*}%
Define the usual Wald test statistic:%
\begin{equation}
W_{T}=Tg^{\prime }(\hat{\theta}_{T})\left[ \frac{\partial g}{\partial \theta
^{\prime }}(\hat{\theta}_{T})\hat{V}_{T}\frac{\partial g^{\prime }}{\partial
\theta }(\hat{\theta}_{T})\right] ^{-1}g(\hat{\theta}_{T}).
\label{wald stat}
\end{equation}

For linear $g$ the statistic converges to a $\chi _{q}^{2}$ distribution,
for other, e.g. polynomial restrictions the limit distribution may be not $%
\chi ^{2}.$ The limit results for the statistic for testing general
polynomial restrictions can be found in (Dufour et al., 2013, 2017); it is
also established there that under some conditions the statistic may diverge
when $q>1$. Below is an example of divergence.

\textbf{Example 1. Restrictions for which the Wald statistic diverges.}

\textit{Consider for }$\theta =\left( \mathbf{x,y,z,w}\right) ^{\prime }$%
\textit{\ the set of restrictions, }$H_{0}:$%
\begin{equation*}
\left\{ 
\begin{array}{ccc}
\mathbf{xy} & = & 0; \\ 
\mathbf{xw} & = & 0; \\ 
\mathbf{yz} & = & 0.%
\end{array}%
\right.
\end{equation*}

\textit{Then the Wald statistic for }$\hat{\theta}_{T}=\left( x,y,z,w\right)
,$\textit{\ assuming that the covariance matrix is identity }$\hat{V}%
=V=I_{4},$\textit{\ is }$\ $%
\begin{equation*}
W=T\left( w^{2}+y^{2}\right) \frac{x^{2}+z^{2}}{w^{2}+x^{2}+y^{2}+z^{2}}.
\end{equation*}

\textit{Suppose that the true parameter value is }$\bar{\theta}=(0,0,1,1);$%
\textit{\ }$H_{0}$\textit{\ then holds. Suppose that the estimated parameter 
}$\hat{\theta}_{T}=\left( x,y,z,w\right) $\textit{\ as }$T\rightarrow \infty 
$\textit{\ is consistent and satisfies }%
\begin{equation*}
z\rightarrow _{p}1;w\rightarrow _{p}1;T^{\frac{1}{2}}x\rightarrow
_{p}Z_{1};T^{\frac{1}{2}}y\rightarrow _{p}Z_{2}
\end{equation*}%
\textit{where }$Z_{1},Z_{2}$\textit{\ are independent standard normals}$.$%
\textit{\ Then the Wald statistic can be expressed as}%
\begin{equation*}
W=\left( T+Z_{2}^{2}+o_{p}\left( 1\right) \right) \frac{Z_{1}^{2}+T+o_{p}%
\left( 1\right) }{2T+Z_{1}^{2}+Z_{2}^{2}+o_{p}\left( 1\right) }.
\end{equation*}

\textit{This is}

\begin{equation*}
\mathit{W=T+O}_{p}\left( 1\right) \mathit{.}
\end{equation*}

\textit{As }$T\rightarrow \infty $ \textit{the statistic diverges under }$%
H_{0}.$

We shall assume that each $g_{l}\left( \theta \right) $\ is a polynomial of
order $m_{l}$\ in the components of $\theta .$ Then for any $\bar{\theta}$
each polynomial component $g_{l}\left( \theta \right) ,$ can be written
around $\bar{\theta}$ as%
\begin{equation}
g_{l}(\theta )=\sum_{\gamma =0}^{m_{l}}\sum_{j_{1}+...+j_{p}=\gamma
}c_{l}(j_{1},...,j_{p},\bar{\theta})\dprod\limits_{k=1}^{p}\left( \theta
_{k}-\bar{\theta}_{k}\right) ^{j_{k}}  \label{power}
\end{equation}%
with some coefficients $c(j_{1},...,j_{p},\bar{\theta}).$

If the value $\bar{\theta}$ satisfies the null hypothesis, then%
\begin{equation*}
c_{l}(0,...,0,\bar{\theta})=0
\end{equation*}%
for each $l.$

A polynomial function is eiher identically zero or non-zero a.e. with
respect to the Lebesgue measure. Consider a square matrix $G(y)$ of
polynomials of variable $y\in \mathbb{R}^{p}.$ We say that the polynomial
matrix $G(y)$\ is non-singular if its determinant is a non-zero polynomial.

The rank of the $q\times p$\ matrix $G(y)$ is the largest dimension of a
square non-singular submatrix.

Unlike matrices of constants for polynomial matrices the rows may be
linearly independent vectors of polynomial functions, while the matrix may
have defficient rank. For example, in the matrix%
\begin{equation*}
\left( 
\begin{array}{cc}
y & 0 \\ 
y^{2} & 0%
\end{array}%
\right)
\end{equation*}%
the two rows are given by independent vectors of polynomials, but the rank
of this matrix of polynomials is one.

\textbf{Assumption 2.} \textit{The }$q\times 1$\textit{\ function }$g\left(
\theta \right) $\textit{\ is a vector of polynomial functions; the matrix of
polynomials }$G\left( \theta \right) =\frac{\partial g}{\partial \theta
^{\prime }}\left( \theta \right) $\textit{\ is of rank }$q.$

This does not exclude the possibility of reduced rank at some particular
point or on a low dimensional space.

Under the stated assumptions for $g\left( \bar{\theta}\right) =0$ the
standard asymptotic $\chi _{q}^{2}$ distribution holds for $W_{T}$ as long
as $G\left( \bar{\theta}\right) =\frac{\partial g}{\partial \theta ^{\prime }%
}(\bar{\theta})$ is a (numerical) matrix of rank $q.$

Each restriction $g_{l}\left( \theta \right) $ can be represented as a sum 
\begin{equation}
g_{l}\left( \theta \right) =\bar{g}_{l}\left( \theta -\bar{\theta}\right)
+r_{l}\left( \theta -\bar{\theta}\right) ,  \label{decomp g}
\end{equation}%
where $\bar{g}_{l}\left( .\right) $ denotes the lowest degree non-zero
homogeneous polynomial and has degree $\bar{\gamma}_{l}+1,$ the degrees of
all non-zero monomials in $r_{l}\left( .\right) $ are \TEXTsymbol{>}$\bar{%
\gamma}_{l}+1.$ We ascribe the degree of homogeneity $\infty $ to a function
that is identically zero.

Correspondingly to $\left( \ref{decomp g}\right) $, in the matrix 
\begin{equation*}
G\left( \theta \right) =\frac{\partial g}{\partial \theta ^{\prime }}(\theta
)
\end{equation*}%
for each row write 
\begin{equation*}
G_{l}\left( \theta \right) =\bar{G}_{l}\left( \theta -\bar{\theta}\right)
+R_{l}\left( \theta -\bar{\theta}\right) ,
\end{equation*}%
the degree of any non-zero homogenious polynomial in the row vector, $\bar{G}%
_{l}\left( \theta -\bar{\theta}\right) ,$ is $\bar{\gamma}_{l}$; any
non-zero monomial in $R_{l}\left( \theta -\bar{\theta}\right) $ has degree
higher than $\gamma _{l}.$ Then collecting the lowest degree homogeneous
polynomials in each row we have 
\begin{equation}
G\left( \theta \right) =\bar{G}\left( \theta -\bar{\theta}\right) +R\left(
\theta -\bar{\theta}\right) .  \label{decomp G}
\end{equation}

\section{\textbf{The property of full rank reached at lowest degrees (FRALD)
and FRALD-T}}

\bigskip \textbf{Definition (FRALD). }\textit{If the matrix }$\bar{G}\left(
\theta -\bar{\theta}\right) $\textit{\ of lowest degree polynomials for }$%
g\left( \theta \right) $\textit{\ is of full rank }$q$\textit{\ we say that
the Full Rank at Lower Degrees (FRALD) property is satisfied for }$g\left(
.\right) $\textit{\ and }$\bar{\theta}.$

Examples in Dufour et al (2017) illustrate the possibilities that the FRALD
property may hold at some points $\bar{\theta},$ but not others, and that
even if FRALD property does not hold for $g\left( .\right) $ at $\bar{\theta}%
,$ it may hold for $Sg\left( .\right) ,$ where $S$ is a non-degenerate
numerical matrix.

Recall that the distribution of the Wald statistic is invariant with respect
to non-degenerate linear transformation of the restrictions.

\textbf{Definition (FRALD-T). }\textit{There exists some numerical
non-degenerate matrix }$S$\textit{\ such that FRALD holds for }$Sg\left(
\theta \right) $\textit{\ at }$\bar{\theta}.$

If FRALD-T holds for $g,$ then for some $S$ FRALD holds for $Sg,$ meaning
that $\overline{SG}\left( .\right) $ is a full rank matrix of polynomials.

It is shown in Dufour et al (2017) that for polynomial $g\left( x\right) $
with a full rank matrix $G\left( x\right) $ there always exists a
non-degenerate numerical matrix $S,$ such that $SG\left( x\right) $ has the
property that $\overline{SG}\left( x\right) $ has all the rows represented
by linearly independent vectors of polynomials (each row contains non-zero
homogeneous polynomials); these rows could be stacked by a permutation in an
"eschelon form", with the degrees of the non-zero homogeneous polynomials in
non-decreasing order.

The eschelon form is given by%
\begin{equation}
\overline{SG}\left( x\right) =\left[ 
\begin{array}{c}
\left[ \overline{SG}\left( x\right) \right] _{1} \\ 
\vdots \\ 
\left[ \overline{SG}\left( x\right) \right] _{i} \\ 
\vdots \\ 
\left[ \overline{SG}\left( x\right) \right] _{v}%
\end{array}%
\right] ,  \label{eschelon}
\end{equation}%
where $\left[ \overline{SG}\left( x\right) \right] _{i}$ has dimension $%
n_{i}\times q,$ all non-zero polynomials in $\left[ \overline{SG}\left(
x\right) \right] _{i}$ have degree $\bar{s}_{i}$ for $i=1,\ldots ,\,\nu $,
with $0\leq \bar{s}_{1}<\cdots <\bar{s}_{i}<\cdots <\bar{s}_{\nu },$\ and
all the rows of $\overline{SG}\left( x\right) $ are linearly independent
functions. Once any $S$ that provides such a structure is found, the rank of 
$\overline{SG}\left( \theta \right) $ is either $q,$ and FRALD-T holds, or
is less than $q,$ in which case this property is violated. An algorithm to
find $S$ is provided in Dufour et al (2017).

\textbf{Example 2 (Example 1 continued). FRALD-T does not hold.}

\textit{Take for }$\theta ^{\prime }=\left( x;y;w;z\right) $\textit{\ the
function }$g\left( \theta \right) =\left( 
\begin{array}{c}
xy \\ 
xw \\ 
yz%
\end{array}%
\right) .$\textit{\ With }$\bar{\theta}^{\prime }=\left( 0,0,1,1\right) $%
\textit{\ denote }$\theta ^{\prime }=\left( x;y;1+\tilde{w};1+\tilde{z}%
\right) $\textit{\ with }$\sqrt{T}\left( x;y;\tilde{w};\tilde{z}\right)
^{\prime }$\textit{\ converging to }$N\left( 0,I_{4}\right) .$\textit{\ Then 
}%
\begin{equation*}
G\left( \theta \right) =\left( 
\begin{array}{cccc}
y & x & 0 & 0 \\ 
1+\tilde{w} & 0 & x & 0 \\ 
0 & 1+\tilde{z} & 0 & y%
\end{array}%
\right) ;
\end{equation*}%
\textit{by applying a transformation (here permutation), }$P,$\textit{\ to
the rows of this matrix we get }%
\begin{eqnarray*}
PG\left( \theta \right) &=&\left( 
\begin{array}{cccc}
1+\tilde{w} & 0 & x & 0 \\ 
0 & 1+\tilde{z} & 0 & y \\ 
y & x & 0 & 0%
\end{array}%
\right) ; \\
&&\text{\textit{with the eschelon form}} \\
\overline{PG}\left( \theta \right) &=&\left( 
\begin{array}{cccc}
1 & 0 & 0 & 0 \\ 
0 & 1 & 0 & 0 \\ 
y & x & 0 & 0%
\end{array}%
\right) .
\end{eqnarray*}%
\textit{The matrix }$\overline{PG}\left( \theta \right) $\textit{\ has
independent polynomial row vectors, and the rows are stacked so that the
degrees of "leading" polynomials do not decline from row to row (eschelon
form). The rank of the matrix }$\overline{PG}\left( \theta \right) $\textit{%
\ is not full in an eschelon form, no linear transformation applied to }$G$ 
\textit{can remedy this rank defficiency. So FRALD-T does not hold for this
example.}

In Dufour et al (2014, 2017) the limit distribution for the Wald statistic
was established for $\bar{\theta}$ when the FRALD-T property holds.

The example 1 here of the case where the statistic was shown to diverge does
not satisfy FRALD-T. The next section demonstrates the mechanism whereby the
violation of the FRALD-T property leads to divergence of the statistic.

\section{Divergence of the Wald statistic where FRALD-T does not hold}

Assume that for $\bar{\theta}$ for which the null is satisfied, $g\left( 
\bar{\theta}\right) =0,$ the FRALD-T property does not hold. Without loss of
generality we may assume that $G\left( x\right) =\frac{\partial g}{\partial
x^{\prime }}$ is such that the eschelon form $\left( \ref{eschelon}\right) $
applies to $\bar{G}\left( x\right) $ (so that $S$ in FRALD-T and in $\left( %
\ref{eschelon}\right) $ is identity).

Denote by $Y$ the Gaussian limit $\sqrt{T}(\hat{\theta}-\bar{\theta}%
)\rightarrow _{d}Y\sim N\left( 0,V\right) .$ Define $\Delta _{T}:=\mathrm{%
diag}[T^{s_{1}/2}I_{n_{1}},\ldots ,\,T^{s_{q}/2}I_{n_{v}}].$ With the
scaling $\Delta _{T}$ we get

\begin{equation}
T^{1/2}\Delta _{T}\,g(\hat{\theta}_{T})\underset{T\rightarrow \infty }{%
\overset{d}{\rightarrow }}\bar{g}\left( Y\right) \,,  \label{Delta}
\end{equation}%
\begin{equation}
\Delta _{T}\,G(\hat{\theta}_{T})\hat{V}_{T}G(\hat{\theta}_{T})^{\prime
}\Delta _{T}\underset{T\rightarrow \infty }{\overset{d}{\rightarrow }}\bar{G}%
\left( Y\right) V\bar{G}\left( Y\right) ^{\prime }\,,  \label{Delta2}
\end{equation}%
where $\bar{G}\left( Y\right) V\bar{G}\left( Y\right) ^{\prime }$ has rank $%
r<q$ when FRALD-T does not hold. Consequently, inverting the consistent
estimator $\left[ \Delta _{T}\,G(\hat{\theta}_{T})\hat{V}_{T}G(\hat{\theta}%
_{T})^{\prime }\Delta _{T}\right] $ for 
\begin{equation*}
W=T^{1/2}\Delta _{T}\,g(\hat{\theta}_{T})^{\prime }\left[ \Delta _{T}\,G(%
\hat{\theta}_{T})\hat{V}_{T}G(\hat{\theta}_{T})^{\prime }\Delta _{T}\right]
^{-1}T^{1/2}\Delta _{T}\,g(\hat{\theta}_{T})
\end{equation*}%
will lead to an explosion as $T\rightarrow \infty .$

We next examine the matrix $\bar{\Sigma}_{T}\left( \hat{\theta}%
_{T},V_{T}\right) =\left[ \Delta _{T}\,G(\hat{\theta}_{T})\hat{V}_{T}G(\hat{%
\theta}_{T})^{\prime }\Delta _{T}\right] $ and its limit eigenvalues which
provide the key ingredient to prove the divergence of the Wald statistic.

Denote by $\bar{\lambda}_{1T}\left( \hat{\theta}_{T}\right) ,\bar{\lambda}%
_{2T}\left( \hat{\theta}_{T}\right) ,...,\bar{\lambda}_{qT}\left( \hat{\theta%
}_{T}\right) $ the eigenvalues of the matrix $\bar{\Sigma}_{T}\left( \hat{%
\theta}_{T},V_{T}\right) ,$ arranged in decreasing order: $\bar{\lambda}%
_{1T}\left( \hat{\theta}_{T}\right) \geq \bar{\lambda}_{2T}\left( \hat{\theta%
}_{T}\right) \geq ...\geq \bar{\lambda}_{qT}\left( \hat{\theta}_{T}\right) $
and denote by 
\begin{equation}
\bar{\Lambda}_{T}=diag\left[ \bar{\lambda}_{1T}\left( \hat{\theta}%
_{T}\right) ,\bar{\lambda}_{2T}\left( \hat{\theta}_{T}\right) ,...,\bar{%
\lambda}_{qT}\left( \hat{\theta}_{T}\right) \right]  \label{diag}
\end{equation}
the $q\times q$ diagonal matrix of these eigenvalues.

We prove several auxilliary results about eigenvalues of non-random
polynomial matrices (proofs are in the next section).

Start with $B\left( x,U\right) =G(x)UG(x)^{\prime },$ where $G\left(
x\right) $ is a non-zero matrix of polynomial functions and define the
characteristic polynomial, $p\left( \lambda ;B\left( x,U\right) \right)
=\det \left[ \lambda I_{p}-B\left( x,U\right) \right] $.

The next proposition describes a polynomial representation for the
coefficients of $p\left( \lambda ;B\left( x,U\right) \right) $ as a
polynomial in $\lambda $. Denote by $F_{p}$ the set of all real symmetric
positive-definite matrices.

\textbf{Proposition 1.} \textit{Let }$x\in R^{p}$\textit{\ and }$G\left(
x\right) $\textit{\ be a }$q\times p$\textit{\ non-zero matrix of polynomial
functions in }$x$\textit{\ such that rank}$\left[ G\left( x\right) \right]
=q $\textit{\ a.e., }$U\in F_{p}$\textit{, }$B\left( x,U\right)
=G(x)UG(x)^{\prime }$\textit{\ and }$p_{B}\left( \lambda ;x,U\right) =\det %
\left[ \lambda I_{p}-B\left( x,U\right) \right] $\textit{\ is the
characteristic polynomial. Then }$p_{B}\left( \lambda ;x,U\right) $\textit{\
can be written as}%
\begin{equation}
p_{B}\left( \lambda ;x,U\right) =\lambda ^{q}+\Sigma _{k=0}^{q}a_{k}\left(
x,U\right) \lambda ^{q-k},  \label{charact polyn}
\end{equation}%
\textit{where the coefficients }$a_{k}\left( x,U\right) $\textit{\ have the
following polynomial expansions}%
\begin{equation}
a_{k}\left( x,U\right) =D_{m_{k}\left( U\right) }\left( x,U\right) +\tilde{R}%
_{k}\left( x,U\right) ,  \label{Dm(k)}
\end{equation}%
\textit{where }$D_{m_{k}\left( U\right) }\left( x,U\right) $\textit{\ is a
homogeneous in }$x$\textit{\ polynomial of degree }$m_{k}\left( U\right) ,$%
\textit{\ and }$\tilde{R}_{k}\left( x,U\right) $\textit{\ is a sum of
polynomials with any non-zero mononomials of degree strictly greater than }$%
m_{k}\left( U\right) .$\textit{\ Further, if }$m_{k}=\underset{Q\in F_{p}}{%
\min }m_{k}\left( Q\right) ,$\textit{\ then }$m_{k}\left( U\right) =m_{k}$%
\textit{\ for almost all }$U\in F_{p}.$

\textbf{Example 3 Restrictions of Example 1 but with a covariance matrix }$U$%
\textbf{\ for which }$m_{k}\left( U\right) >m_{k}.$

\textit{In the example 1 we had divergence at the rate }$T$\textit{\ when
the matrix }$U$\textit{\ was identity. The characteristic polynomial for the
same restrictions, }$\mathit{g}\left( .\right) $\textit{\ of example 1 with
a covarince matrix }$U,$\textit{\ possibly different from }$I,$ \textit{has
as the }$q-th$\textit{\ (here for }$q=3)\,$\textit{\ coefficient the
determinant of }$G(x)UG(x)^{\prime }.$\textit{\ With }$U=I$\textit{\ the
determinant of }$G(x)IG(x)^{\prime }$\textit{\ provides }%
\begin{equation*}
a_{3}\left( x,I\right)
=w^{2}x^{2}y^{2}+2wx^{2}y^{2}+x^{4}y^{2}+x^{2}y^{4}+x^{2}y^{2}z^{2}+%
\allowbreak 2x^{2}y^{2}z+2x^{2}y^{2},
\end{equation*}%
\textit{and we note that the lowest degree monomial is }$x^{2}y^{2}.$\textit{%
\ It can be verified that for these restrictions and }$\bar{\theta}$\textit{%
\ we get }$m_{3}=4$\textit{\ (so that there can be no }$U$\textit{\ for
which the degree could be smaller) and by Proposition 1 }$m_{3}\left(
U\right) =m_{3}=4$\textit{\ \ for almost every }$U\in F_{q}.$ \textit{%
However, below we provide }$U$\textit{\ for which }$a_{3}\left( x,U\right) $
is such that $m_{3}\left( U\right) >4.$\textit{\ Consider}%
\begin{equation*}
U=\left( 
\begin{array}{cccc}
1 & \sqrt{.98} & 0 & 0 \\ 
\sqrt{.98} & 1 & .1 & .1 \\ 
0 & .1 & 1 & 0 \\ 
0 & .1 & 0 & 1%
\end{array}%
\right) .
\end{equation*}%
\textit{For this }$U$\textit{\ we get }%
\begin{eqnarray*}
a_{3}\left( x,U\right) &=&0.01w^{2}x^{2}y^{2}-0.197\,99\allowbreak
wx^{3}y^{2}-0.2wx^{2}y^{3}-0.02wx^{2}y^{2}z+0.98x^{4}y^{2} \\
&&+1.\,\allowbreak 979\,9\allowbreak x^{3}y^{3}+0.197\,99x^{3}\allowbreak
y^{2}z+x^{2}y^{4}+0.2x^{2}y^{3}z+0.01x^{2}y^{2}z^{2}
\end{eqnarray*}%
\textit{with the lowest degree of monomial }$m_{3}\left( U\right)
=6>m_{3}=4. $

Since the coefficients of a polynomial represent symmetric polynomials in
the roots (e.g., Horn and Johnson, 1985, Section 1.2), elementary symmetric
polynomials in the eigenvalues can be expressed as polynomial functions in $%
x.$ Denote by $I_{q}\left( k\right) $ the set of all combinations of $k$
integers out of $\left\{ 1,...,q\right\} .$ Denote by $P_{k}\left( \lambda
_{1},...\lambda _{q}\right) $ the $k-th$ elementary symmetric polynomial in $%
\lambda _{1},...\lambda _{q}:$%
\begin{equation*}
P_{k}\left( \lambda _{1},...\lambda _{q}\right) =\Sigma _{\left\{
i_{1},...,i_{k}\right\} \in I_{q}\left( k\right) }\Pi _{j=1}^{k}\lambda
_{i_{j}}.
\end{equation*}

\textbf{Corollary to Proposition 1. }\textit{For the eigenvalues }$\lambda
_{i}\left( x,U\right) ,$\textit{\ }$i=1,...,q,$\textit{\ \ that are the
solutions of the characteristic polynomial, we have that }%
\begin{equation}
P_{k}\left[ \lambda _{1}\left( x,U\right) ,...\lambda _{q}\left( x,U\right) %
\right] =\left( -1\right) ^{k}a_{k}\left( x,U\right)  \label{sym polyn}
\end{equation}%
\textit{and thus the representation }$\left( \ref{Dm(k)}\right) $\textit{\
applies. }

In the next proposition we apply scaling to the argument $x$ by considering $%
x=T^{-1/2}y$ and exploit the polynomial terms from $\left( \ref{Dm(k)}%
\right) $ with lowest degree of homogeneity in $\left( \ref{sym polyn}%
\right) $ to establish the rates for the eigenvalues of a scaled polynomial
matrix. Recall that from the convergence result $\left( \ref{Delta2}\right) $
the matrix scaling $\Delta _{T}$ is associated with $G\left( .\right) .$ We
show that when rank of $\bar{G}\left( .\right) $ is less than $q$ (in
violation of the FRALD-T condition) some eigenvalues will be converging to
zero and additional scaling can be applied to have the eigenvalues converge
to continuous limit functions. This additional scaling will provide the
divergence rate.

\textbf{Proposition 2.} \textit{Under the conditions of Proposition 1
consider the scaled matrix for }$y\in R^{p}:$\textit{\ }%
\begin{equation*}
M_{T}\left( y,U\right) =\Delta _{T}G\left( T^{-1/2}y\right) UG\left(
T^{-1/2}y\right) ^{\prime }\Delta _{T}
\end{equation*}%
\textit{and its eigenvalues }$\bar{\lambda}_{l}^{\left( T\right) }\left(
y,U\right) =\bar{\lambda}_{l}^{\left( T\right) }\left( T^{-1/2}y,U\right) ,$%
\textit{\ }$l=1,...,q$\textit{\ in descending order. Then for some
non-negative integers }$\beta _{l}=\beta _{l}\left( U\right) $\textit{\ that
satisfy}%
\begin{eqnarray*}
\beta _{l} &=&0\text{ for }1\leq l\leq r; \\
\beta _{l} &\geq &1\text{ for }l>r,
\end{eqnarray*}%
\textit{we have that}%
\begin{equation*}
T^{\beta _{l}}\bar{\lambda}_{l}^{\left( T\right) }\left( y,U\right)
\rightarrow \lambda _{l}\left( y,U\right) \text{ for almost all }y\in
R^{p},l=1,...,q,
\end{equation*}%
\textit{where }$\lambda _{l}\left( y,U\right) $\textit{\ are continuous a.e.
non-zero functions.}

Thus we see that for eigenvalues beyond $r$ additional non-trivial scaling
provides convergence to a continuous a.e. non-zero function.

The next proposition shows that convergence with these rates to a continuous
(but now in some exceptional cases possibly zero) function is preserved when 
$U$ is replaced with a sequence $U_{T},$ of matrices from $F_{p}$ such that $%
U_{T}\rightarrow U.$

\textbf{Proposition 3.} \textit{Under the conditions of Proposition 2
consider a sequence }$U_{T}\in F_{p},$\textit{\ such that }$U_{T}\rightarrow
U$\textit{. Then}

\textit{(a) if }$m_{l}\left( U\right) =m_{l}$\textit{\ for }$l=1.,,,q$%
\textit{\ we have that}%
\begin{eqnarray*}
T^{\beta _{l}}\bar{\lambda}_{l}^{\left( T\right) }\left( y,U_{T}\right)
&\rightarrow &\lambda _{l}\left( y,U\right) \text{ for almost all }y\in
R^{p},l=1,...,q; \\
\lambda _{l}\left( y,U\right) &>&0\text{ a.e.;}
\end{eqnarray*}

\textit{(b) if for some }$k\in \left\{ 1,...,q\right\} $\textit{\ it holds
that }$m_{l}\left( U\right) =m_{l}$\textit{\ for }$l<k$\textit{\ and }$%
m_{k}<m_{k}\left( U\right) ,$\textit{\ then }%
\begin{equation*}
T^{\beta _{l}}\bar{\lambda}_{l}^{\left( T\right) }\left( y,U_{T}\right)
\rightarrow 0\text{ with }\beta _{l}=\frac{1}{2}\left( m_{k}-m_{k-1}\right)
\geq 1\text{ for }l\geq k.
\end{equation*}%
Recall that case (a) will hold for almost all $U$ by Proposition 1. Example
3 illustrates part (b) of Proposition 2: there $m_{3}=4,$ so $\beta _{3}=2,$
but $m_{3}\left( U\right) =6>4$ and thus $\lambda _{3}^{\left( T\right)
}\left( T^{-1/2}y,U_{T}\right) $ scaled up by $T^{2}$ goes to zero. It is
possible that for \thinspace $U_{T}$ itself $m_{3}\left( U_{T}\right) =6;$
in that case $T^{3}\lambda _{3}^{\left( T\right) }\left(
T^{-1/2}y,U_{T}\right) $ converges to a non-zero limit. Alternatively, if $%
U_{T}$ has $m_{3}\left( U_{T}\right) =m_{3}$ for almost all $U_{T}$ but $%
U_{T}$ converges to $U$ sufficiently fast the rate could still be as high as 
$T^{3}.$ However, in case (b) to get a precise rate we also need to consider
the convergence rate for $U_{T}.$

The next proposition applies the deterministic properties to provide limits
for eigenvalues of the random matrix $\bar{\Sigma}_{T}\left( \hat{\theta}%
_{T},\hat{V}_{T}\right) ,$ in the diagonal eigenvalue matrix $\bar{\Lambda}%
_{T}\left( \hat{\theta}_{T}\right) $ of $\left( \ref{diag}\right) $. Without
loss of generality consider $\bar{\theta}=0.$

\textbf{Proposition 4.} \textit{Suppose that assumptions 1,2 hold at }$\bar{%
\theta}=0.$\textit{\ Then there is a sequence of integers }$\beta
_{l},l=1,...,q,$\textit{\ which depends on }$G$\textit{\ and }$V,$\textit{\
such that}%
\begin{eqnarray*}
\beta _{l} &=&0\text{ for }l=1,..,r, \\
\bar{\beta} &=&\underset{r<l\leq q}{\max }\beta _{l}\geq 1; \\
\tilde{\Delta}_{T} &=&diag\left[ T^{\beta _{1}/2},...,T^{\beta _{q}/2}\right]
; \\
&&\tilde{\Delta}_{T}\bar{\Lambda}_{T}\left( \hat{\theta}_{T}\right) \tilde{%
\Delta}_{T}\underset{d}{\rightarrow }diag\left[ \lambda _{1}\left( Y\right)
,...,\lambda _{T}\left( Y\right) \right]
\end{eqnarray*}%
\textit{with all }$\lambda _{l}\left( y\right) $\textit{\ continuous
non-negative functions a.e..}

We see that if FRALD-T were not violated, no additional scaling would be
required, but once it is violated the extra scaling is captured by $\bar{%
\beta}_{l}\geq 1$ for $r<l\leq q$ that determines the rate of explosion of
the Wald statistic. The Theorem below shows this.

\textbf{Theorem.} \textit{Under the conditions of Proposition 4 if FRALD-T
property does not hold, i.e. }$r<q,$\textit{\ then we have for }$\bar{\beta}%
\geq 1$\textit{\ that }%
\begin{equation*}
W_{T}=W_{T}\left( \hat{\theta}_{T};g,V_{T}\right) >T^{\bar{\beta}}\mu
_{T}\left( \hat{\theta}_{T},V_{T}\right) ,
\end{equation*}%
\textit{where }$\mu _{T}\left( \hat{\theta}_{T},V_{T}\right) \underset{d}{%
\rightarrow }\mu \left( Y\right) ,$\textit{\ a continuous positive a.e.
function.}

We thus see that if FRALD-T is violated the rate of the exlosion is at least 
$T$ (as in example 1 here), but could be stronger even with the same
restrictions (as could be in example 3).

\section{Proofs}

\textbf{Proof of Proposition 1.}

First, consider the polynomial expansion for $p_{B}\left( \lambda
;x,U\right) $, given e.g. in Harville (2008, Corollary 13,7,4). Denote by $%
I_{q}\left( k\right) $ the set of all combinations of $k$ integers out of $%
\left\{ 1,...,q\right\} ;$ denote for $\left\{ i_{1},...,i_{r}\right\} \in
I_{q}\left( r\right) $ by $B^{\left\{ i_{1},...,i_{r}\right\} }$ a minor of
the matrix $B$ obtained by striking out all the rows and columns numbered $%
i_{1},...,i_{r}.$ Then 
\begin{eqnarray*}
P_{B}\left( \lambda ;x,U\right)  &=&\Sigma _{r=0}^{q}\left( -1\right)
^{q-r}\lambda ^{r}\Sigma _{\left\{ i_{1},...,i_{r}\right\} \in I_{q}\left(
r\right) }\det \left[ B\left( x,U\right) ^{\left\{ i_{1},...,i_{r}\right\} }%
\right] ; \\
a_{k}\left( x,U\right)  &=&\left( -1\right) ^{q-k}\Sigma _{\left\{
i_{1},...,i_{k}\right\} \in I_{q}\left( k\right) }\det \left[ B\left(
x,U\right) ^{\left\{ i_{1},...,i_{k}\right\} }\right] .
\end{eqnarray*}%
Since all the components of the $B\left( x,U\right) $ matrix are polynomials
in $x$ it follows that the determinants of the minors are also polynomials
in $x.$ Then, given $U,$ denote by $D_{m_{k}\left( U\right) }\left(
x,U\right) $ the homogeneous polynomial in $a_{k}\left( x,U\right) $ of the
lowest degree, denoted $m_{k}\left( U\right) ,$ to obtain the polynomial
expansion of the Proposition.

Next, note that $a_{k}\left( x,U\right) $ is also a polynomial function in
the components of the matrix $U.$ By varying $U$ over $F_{q}$ we can find
the minimum possible $m_{k}\left( U\right) ,$ denoted $m_{k}.$ Thus there is
some matrix, $Q\in F_{q}$ such that for $a_{k}\left( x,Q\right) $ we have
that $m_{k}\left( Q\right) =m_{k}.$ This implies that in $a_{k}\left(
x,Q\right) $ there is a homogeneous polynomial of degree $m_{k}$ in $x$ that
is non-zero, thus has at least one non-zero coefficient on a monomial term
of degree $m_{k}.$ Since this coefficient is a polynomial function of the
components of $Q,$ considering this polynomial over the corresponding
components of all $U\in F_{q}$ we note that it is non-zero a.e.. This
implies that $m_{k}\left( U\right) =m_{k}$ for almost all $U\in F_{q}.$

$\blacksquare $

\textbf{Proof of Proposition 2.}

For every $T$ consider $T^{-1/2}y$ in place of $x$ and $M_{T}\left(
y,U\right) $ in place of $B\left( x,U\right) $ in Proposition 1. Then for
the corresponding characteristic polynomial, $p_{M_{T}}\left( \lambda
;y,U\right) ,$ the expansion similar to $\left( \ref{charact polyn}\right) $
will provide coeffficients%
\begin{equation*}
\tilde{a}_{k}\left( y,U\right) =\left( -1\right) ^{q-k}\Sigma _{\left\{
i_{1},...,i_{k}\right\} \in I_{q}\left( k\right) }\det \left[ M_{T}\left(
y,U\right) ^{\left\{ i_{1},...,i_{k}\right\} }\right] .
\end{equation*}%
Note that we have that $\det \left[ M_{T}\left( y,U\right) ^{\left\{
i_{1},...,i_{k}\right\} }\right] \rightarrow \det \left[ \bar{G}\left(
y\right) V\bar{G}\left( y\right) ^{\prime \left\{ i_{1},...,i_{k}\right\} }%
\right] ,$ that is a non-zero constant for $k\leq r$ and zero for $k>r.$
Therefore the coefficients can be represented as 
\begin{equation*}
\tilde{a}_{k}\left( y,U;T\right) =D_{\tilde{m}_{k}\left( U\right) }\left(
y,U;T\right) +\tilde{R}_{k}\left( y,U;T\right) ,
\end{equation*}%
where $\tilde{m}_{k}\left( U\right) $ is zero for $k=1,...,r.$ But for $%
k=r+1,...,q$ there is some $\gamma _{k}\geq 1$ such that $D_{\tilde{m}%
_{k}\left( U\right) }\left( y,U;T\right) =T^{-\gamma _{k}}\bar{R}_{k}\left(
y\right) ,$ where $\gamma _{k}=\frac{1}{2}\tilde{m}_{k}\left( U\right) $ and 
$\bar{R}_{k}\left( y\right) $ is a (positive a.e.) homogeneous polynomial in 
$y$ of degree $\gamma _{k}.$ So altogether we can write

\begin{equation*}
\tilde{a}_{k}\left( y,U;T\right) =T^{-\gamma _{k}}\bar{R}_{k}\left( y\right)
+\tilde{R}_{k}\left( y,U;T\right) ,
\end{equation*}%
where $\tilde{R}_{k}\left( .\right) $ is a polynomial that can contain
non-zero monomials only of degree strictly higher than $\gamma _{k}.$ Then
apply the representation $\left( \ref{sym polyn}\right) $ to the
corresponding coefficients to write for every $k=1,...,q$ 
\begin{eqnarray*}
P_{k}\left[ \lambda ^{\left( T\right) }(T^{-1/2}y)\right] &=&T^{-\gamma _{k}}%
\bar{R}_{k}\left( y\right) +\tilde{R}_{k}\left( T^{-1/2}y\right) ; \\
\gamma _{k} &=&0,\text{ if }k=1,...,r; \\
\gamma _{k} &\geq &1\text{ for }k=r+1,...,q.
\end{eqnarray*}

The proof is by induction on $k.$

For $k=1$ consider the largest eigenvalue $\lambda _{1}^{\left( T\right)
}(T^{-1/2}y).$ Note that $r\geq 1,$ so that $\gamma _{1}$ is always zero.
Since $P_{1}\left[ .\right] $ is the sum of all eigenvalues we have by
replacing all the $q$ eigenvalues by the largest, $\lambda _{1}^{\left(
T\right) }(T^{-1/2}y),$ that%
\begin{equation*}
q\lambda _{1}^{\left( T\right) }(T^{-1/2}y)\geq \bar{R}_{1}\left( y\right)
+O\left( T^{-1/2}\right) ,
\end{equation*}%
then since $\bar{R}_{1}\left( y\right) >0$ a.e. the limit of $T^{\beta
_{1}}\lambda _{1}^{\left( T\right) }(T^{-1/2}y)$ (with $\beta _{1}=0)$ is
positive a.e..

Suppose that for $k^{\prime }\geq 1$ all $T^{\beta _{l}}\lambda _{l}^{\left(
T\right) }(T^{-1/2}y)$ for $l\leq k^{\prime }$ converge to continuous
positive a.e. functions.

Then by replacing in the symmetric polynomial $P_{k^{\prime }+1}\left[
\lambda ^{\left( T\right) }(T^{-1/2}y)\right] $ all the terms by the
largest, $\lambda _{k^{\prime }+1}^{\left( T\right) }(T^{-1/2}y),$ and
multipying by the rate, $T^{\gamma _{k^{\prime }+1}},$ we can write that%
\begin{eqnarray*}
&&T^{\gamma _{k^{\prime }+1}}\frac{q!}{\left( k^{\prime }+1\right) !\left(
q-k^{\prime }-1\right) !}\left[ \Pi _{l\leq k^{\prime }}\lambda _{l}^{\left(
T\right) }(T^{-1/2}y)\right] \left[ \lambda _{k^{\prime }+1}^{\left(
T\right) }(T^{-1/2}y)\right] \\
&=&\frac{q!}{\left( k^{\prime }+1\right) !\left( q-k^{\prime }-1\right) !}%
\left[ \Pi _{l\leq k^{\prime }}T^{\beta _{l}}\lambda _{l}^{\left( T\right)
}(T^{-1/2}y)\right] \left[ T^{\left( \gamma _{k^{\prime }+1}-\Sigma _{l\leq
k^{\prime }}\beta _{l}\right) }\lambda _{k^{\prime }+1}^{\left( T\right)
}(T^{-1/2}y)\right] \\
&\geq &T^{\gamma _{k^{\prime }+1}}P_{k^{\prime }+1}\left[ \lambda ^{\left(
T\right) }(T^{-1/2}y)\right] \\
&=&\bar{R}_{k^{\prime }+1}\left( y\right) +O\left( T^{-1/2}\right) .
\end{eqnarray*}%
Since the expression in the last line has a limit that is non-zero a.e., so
does the expression in the second line; by the induction hypothesis $\left[
\Pi _{l\leq k^{\prime }}T^{\beta _{l}}\lambda _{l}^{\left( T\right)
}(T^{-1/2}y)\right] $ converges to a continuous positive a.e. function. Thus
for $\beta _{k^{\prime }+1}=\left( \gamma _{k^{\prime }+1}-\Sigma _{l\leq
k^{\prime }}\beta _{l}\right) $ the function $T^{\beta _{k^{\prime
}+1}}\lambda _{k^{\prime }+1}^{\left( T\right) }(T^{-1/2}y)$ converges to a
continuous positive a.e. function.

From the derivation it follows that $\beta _{l}=\frac{1}{2}\tilde{m}%
_{l}\left( U\right) =0$ for $l=1,...,r;$ $\beta _{r+1}=\gamma _{r+1}=\frac{1%
}{2}\tilde{m}_{r+1}\left( U\right) \geq 1;$ and generally $\beta _{k}=\gamma
_{k}-\gamma _{k-1}=\frac{1}{2}\left( \tilde{m}_{k}\left( U\right) -\tilde{m}%
_{k-1}\left( U\right) \right) \geq 1$ for $k>r.\blacksquare $

\textbf{Proof of Proposition 3.}

Under the condition in (a) the degrees of homogeneity $m_{k}\left( U\right) $
for the coefficients of the characteristic polynomial of $M_{T}\left(
y,G,U\right) $ and $m_{k}\left( U_{T}\right) $ for the corresponding
coefficient in $M_{T}\left( y,G,U_{T}\right) $ has to be the same for large
enough $T$ and thus by the proof of Proposition 2 we conclude that $T^{\beta
_{l}}\bar{\lambda}_{l}^{\left( T\right) }\left( y,U_{T}\right) $ have the
same positive a.e. limit as $T^{\beta _{l}}\bar{\lambda}_{l}^{\left(
T\right) }\left( y,U\right) .$

Under the condition in (b) we can write%
\begin{equation*}
T^{m_{k+1}}P_{k+1}\left[ \lambda ^{\left( T\right) }\left(
T^{-1/2}y,U_{T}\right) \right] \geq \left[ \Pi _{i\leq k}T^{\beta
_{l}}\lambda _{l}^{\left( T\right) }\left( T^{-1/2}y,U_{T}\right) \right] %
\left[ T^{(m_{k+1}-m_{k})}\lambda _{k+1}^{\left( T\right) }\left(
T^{-1/2}y,U_{T}\right) \right] ,
\end{equation*}%
where $\left[ \Pi _{i\leq k}T^{\beta _{l}}\lambda _{l}^{\left( T\right)
}\left( T^{-1/2}y,U_{T}\right) \right] $ converges to a function that is
positve a.e., but by the condition for $k$ the left-hand side converges to
zero. Thus $\left[ T^{(m_{k+1}-m_{k})}\lambda _{k+1}^{\left( T\right)
}\left( T^{-1/2}y,U_{T}\right) \right] $ converges to zero and so does $%
\left[ T^{(m_{k+1}-m_{k})}\lambda _{l}^{\left( T\right) }\left(
T^{-1/2}y,U_{T}\right) \right] $ for any $l>k.\blacksquare $

\textbf{Proof of Proposition 4.}

Consider the scaling matrix $\tilde{\Delta}$ and the scaled matrix $\hat{%
\Sigma}_{T}\left( \hat{\theta}_{T},\hat{V}_{T}\right) =\tilde{\Delta}\Delta
_{T}G\left( \hat{\theta}_{T}\right) \hat{V}_{T}G\left( \hat{\theta}%
_{T}\right) ^{\prime }\Delta _{T}\tilde{\Delta}.$ By Assumption 1 since $T^{%
\frac{1}{2}}\left( \hat{\theta}_{T}\right) \underset{d}{\rightarrow }Y$ that
is absolutely continuous by Proposition 2 the eigenvalues of $\tilde{\Delta}%
\Delta _{T}G\left( \hat{\theta}_{T}\right) VG\left( \hat{\theta}_{T}\right)
^{\prime }\Delta _{T}\tilde{\Delta}$ converge in distribution to continuous
functions in $Y,$ some of which are non-zero a.e.. Additionally, for any
sequence of $\hat{V}_{T}$ we can select a subsequence $\hat{V}_{T^{\prime }}$
that converges a.s. to $V$ and by Proposition 3 we have the same result for
the limits of eigenvalues of $\tilde{\Delta}\Delta _{T}G\left( \hat{\theta}%
_{T}\right) \hat{V}_{T}G\left( \hat{\theta}_{T}\right) ^{\prime }\Delta _{T}%
\tilde{\Delta}.\blacksquare $

\textbf{Proof of the Theorem.}

Consider now the matrix $\hat{\Sigma}_{T}\left( \hat{\theta}_{T},\hat{V}%
_{T}\right) =\tilde{\Delta}\Delta _{T}G\left( \hat{\theta}_{T}\right) \hat{V}%
_{T}G\left( \hat{\theta}_{T}\right) ^{\prime }\Delta _{T}\tilde{\Delta}$ as
defined in Proposition 4. The eigenvalues of the scaled matrix, $\hat{\Sigma}%
_{T}\left( \hat{\theta}_{T},\hat{V}_{T}\right) ,$ denoted $\tilde{\lambda}%
_{iT},$ $i=1,...,q$ by Proposition 4 converge in distribution%
\begin{equation*}
\tilde{\lambda}_{i,T}\underset{d}{\rightarrow }\lambda _{i}\left( Y\right) ,%
\text{ }i=1,...,q.
\end{equation*}%
Then rewrite the Wald statistic as%
\begin{equation*}
W_{T}=\left[ \tilde{\Delta}_{T}T^{1/2}\Delta _{T}g\left( \hat{\theta}%
_{T}\right) \right] ^{\prime }\left[ \hat{\Sigma}_{T}\left( \hat{\theta}_{T},%
\hat{V}_{T}\right) \right] ^{-1}\left[ \tilde{\Delta}_{T}T^{1/2}\Delta
_{T}g\left( \hat{\theta}_{T}\right) \right] .
\end{equation*}%
Since $\tilde{\lambda}_{1T}^{-1},...,\tilde{\lambda}_{qT}^{-1}$ are the
eigenvalues of the non-negative definite matrix $\left[ \hat{\Sigma}%
_{T}\left( \hat{\theta}_{T},\hat{V}_{T}\right) \right] ^{-1},$ for any
vector $\xi $ we have%
\begin{equation*}
\xi ^{\prime }\left[ \hat{\Sigma}_{T}\left( \hat{\theta}_{T},\hat{V}%
_{T}\right) \right] ^{-1}\xi \geq \xi ^{\prime }\xi \underset{1\leq i\leq q}{%
\min }\left\{ \tilde{\lambda}_{iT}^{-1}\right\} ,
\end{equation*}%
thus%
\begin{equation*}
W_{T}\geq \frac{\left\Vert \tilde{\Delta}_{T}T^{1/2}\Delta _{T}g\left( \hat{%
\theta}_{T}\right) \right\Vert }{\underset{1\leq i\leq q}{\max }\left\{ 
\tilde{\lambda}_{iT}\right\} }.
\end{equation*}%
We have that 
\begin{equation*}
T^{1/2}\Delta _{T}g\left( \hat{\theta}_{T}\right) \underset{d}{\rightarrow }%
\bar{g}\left( Y\right)
\end{equation*}%
with all components of the vector function $\bar{g}\left( Y\right) $
non-zero a.e. for absolutely continuous $Y.$ Then 
\begin{eqnarray*}
\left\Vert \tilde{\Delta}_{T}T^{1/2}\Delta _{T}g\left( \hat{\theta}%
_{T}\right) \right\Vert ^{2} &=&\Sigma _{i=1}^{r}\left[ T^{\left(
s_{k_{i}}+1\right) /2}g_{i}\left( \hat{\theta}_{T}\right) \right]
^{2}+\Sigma _{i=r+1}^{q}T^{\beta _{i}}\left[ T^{\left( s_{k_{i}}+1\right)
/2}g_{i}\left( \hat{\theta}_{T}\right) \right] ^{2} \\
&\geq &T^{\bar{\beta}}\underset{1\leq i\leq q}{\min }\left\{ \left[
T^{\left( s_{k_{i}}+1\right) /2}g_{i}\left( \hat{\theta}_{T}\right) \right]
^{2}\right\} .
\end{eqnarray*}

Define 
\begin{equation*}
\mu _{T}\left( \hat{\theta}_{T},\hat{V}_{T}\right) =\frac{T^{\bar{\beta}}%
\underset{1\leq i\leq q}{\min }\left\{ \left[ T^{\left( s_{k_{i}}+1\right)
/2}g_{i}\left( \hat{\theta}_{T}\right) \right] ^{2}\right\} }{\underset{%
1\leq i\leq q}{\max }\left\{ \tilde{\lambda}_{iT}\right\} },
\end{equation*}%
then 
\begin{equation*}
W_{T}\geq T^{\bar{\beta}}\mu _{T}\left( \hat{\theta}_{T},\hat{V}_{T}\right) .
\end{equation*}

By Proposition 4, continuity of the eigenvalue function, and of the maximum
of continuous functions 
\begin{equation*}
\underset{1\leq i\leq q}{\max }\left\{ \tilde{\lambda}_{iT}\right\} \underset%
{d}{\rightarrow }\lambda _{\max }\left( Y\right) ,
\end{equation*}%
with the limit functions non-zero a.e.. Also, 
\begin{equation*}
\underset{1\leq i\leq q}{\min }\left\{ \left[ T^{\left( s_{k_{i}}+1\right)
/2}g_{i}\left( \hat{\theta}_{T}\right) \right] ^{2}\right\} \underset{d}{%
\rightarrow }\bar{g}_{i\min }\left( Y\right) ,
\end{equation*}%
which is a piece-wise polynomial continuous function. The ratio 
\begin{equation*}
\mu \left( Y\right) =\frac{\bar{g}_{i\min }\left( Y\right) }{\lambda _{\max
}\left( Y\right) }
\end{equation*}%
exists and is non-zero a.e. and 
\begin{equation*}
\mu _{T}\left( \hat{\theta}_{T},\hat{V}_{T}\right) \underset{d}{\rightarrow }%
\mu \left( Y\right) .
\end{equation*}%
When the FRALD-T condition is violated $\bar{\beta}\geq 1$ and the Wald
statistic diverges to $+\infty .$

$\blacksquare $


\begin{thebibliography}{9}
\bibitem{2013} Dufour, J.-M., Renault, E. and V. Zinde-Walsh, 2013, Wald
tests when restrictions are locally singular, working paper, ArXiV

\bibitem{2017} Dufour, J.-M., Renault, E. and V. Zinde-Walsh, 2017, Wald
tests when restrictions are locally singular, working paper,
https://monde.cirano.qc.ca/\symbol{126}dufourj/Web\_Site/Dufour\_Renault%
\_ZindeWalsh\_2012\_WaldTestsLocallySingularRestrictions\_W.pdf

\bibitem{harvil} Harville, D.A., 2008, Matrix Algebra from a Statistician's
Perspective, Springer-Verlag, New York

\bibitem{horn} Horn, R. G. and Johnson, C. A. (1985), Matrix Analysis,
Cambridge University Press, Cambridge, U.K.
\end{thebibliography}
\end{document}